\documentclass[11pt,a4paper]{amsart}
\usepackage[T1]{fontenc}
\usepackage[utf8]{inputenc}
\usepackage{lmodern}
\usepackage{mathtools,amssymb}
\usepackage{microtype}
\usepackage{enumitem}
\usepackage[left=25mm,right=25mm,top=27mm,bottom=28mm]{geometry}
\usepackage{tikz}
\usetikzlibrary{arrows.meta,positioning}
\usepackage[hidelinks]{hyperref}
\hypersetup{pdftitle={Classification of Isometric Immersions of Spherical and HCMU Footballs into Euclidean Three-Space},pdfauthor={Zhiqiang Wei, Mengyan Xing, Xingya Zhao}}
\numberwithin{equation}{section}
\newtheorem{theorem}{Theorem}[section]
\newtheorem{lemma}[theorem]{Lemma}
\newtheorem{proposition}[theorem]{Proposition}
\newtheorem{corollary}[theorem]{Corollary}
\theoremstyle{definition}
\newtheorem{definition}[theorem]{Definition}
\newtheorem{example}[theorem]{Example}
\theoremstyle{remark}
\newtheorem{remark}[theorem]{Remark}
\newcommand{\R}{\mathbb{R}}
\newcommand{\C}{\mathbb{C}}
\newcommand{\N}{\mathbb{N}}
\newcommand{\Z}{\mathbb{Z}}
\newcommand{\E}{\mathbb{E}}
\newcommand{\dd}{\,\mathrm{d}}
\newcommand{\II}{\mathrm{II}}
\newcommand{\I}{\mathrm{I}}
\newcommand{\tr}{\operatorname{tr}}
\newcommand{\Area}{\operatorname{Area}}
\setlist[enumerate,1]{label=\textup{(\roman*)},leftmargin=2.1em,itemsep=3pt,topsep=5pt}
\setlist[enumerate,2]{label=\textup{(\alph*)},leftmargin=2em,itemsep=2pt,topsep=3pt}
\allowdisplaybreaks[1]
\title[Classification of football immersions]{Classification of Isometric Immersions of Spherical and HCMU Footballs into $\E^{3}$}
\author[Z. Wei]{Zhiqiang Wei}
\address{School of Mathematics and Statistics, Henan University, Kaifeng, China\newline
Center for Applied Mathematics of Henan Province, Henan University, Zhengzhou, China}
\email{10100123@vip.henu.edu.cn}
\date{}
\begin{document}
\begin{abstract}
We study global isometric immersions of rotational metrics with two conical ends into Euclidean three-space. For a metric $g=\mathrm{d}u^2+f(u)^2\mathrm{d}\theta^2$ of positive Gauss curvature, we prove that vanishing endpoint traces of the weighted mixed shape coefficient $f^2b$ force every smooth isometric immersion to be rotational. The resulting congruence classes are parametrized by an integer angular covering degree bounded below by the two cone-angle parameters. The proof combines a scalar elliptic equation derived from the Gauss--Codazzi system, the maximum principle, and a global angular closing condition. For spherical footballs, isolated-singularity theory and a direct estimate in coordinates conformal for the second fundamental form establish the trace condition for every smooth immersion of the regular part. We also classify all rotational isometric immersions of HCMU footballs, including those whose curvature changes sign. When the minimum curvature is nonnegative, this classification exhausts the immersions satisfying the trace condition, and hence all immersions admitting a $C^2$ polar extension at both tips. 
\end{abstract}
\maketitle
\enlargethispage{3pt}

\section{Introduction}
Rotational conical metrics on the sphere provide a natural setting in which to study the relationship between intrinsic symmetry and the geometry of isometric immersions. On the twice-punctured sphere, such a metric takes the form $g=\mathrm{d}u^2+f(u)^2\mathrm{d}\theta^2$, where $f$ vanishes linearly at the endpoints of the meridian interval. The constant-curvature case leads to spherical footballs \cite{Troyanov1989,WeiWuXu2026}, while the nonconstant-curvature case with umbilical curvature Hessian leads to HCMU footballs \cite{ChenChenWu2005}. We study global isometric immersions of these metrics into $\E^3$, including maps that cover their images more than once.

The rotational construction itself has a classical origin. Engman \cite{Engman2004} studies isometric embeddings of $S^1$-invariant metrics on $S^2$ and records the usual profile formula in meridian arclength coordinates. Giblin and Hwang \cite{GiblinHwang2004} give a related construction with an angular parameter, together with a discussion of endpoint geometry and periodicity. Section~\ref{subsec:classical} explains its relation to the formulas used here. These constructions motivate the question of whether the intrinsic circle symmetry must be preserved by an arbitrary isometric immersion.

Spherical footballs belong to the broader class of spherical conical metrics, namely metrics of constant Gauss curvature one with finitely many conical singularities. Troyanov classified the two-cone case \cite{Troyanov1989} and established existence results for prescribed curvature on compact conical surfaces \cite{Troyanov1991}. Luo and Tian treated spherical metrics with cone angles below $2\pi$ \cite{LuoTian1992}. For three singularities, existence criteria were obtained by Umehara and Yamada \cite{UmeharaYamada2000} and by Eremenko \cite{Eremenko2004}. Developing maps and reducible monodromy provide another description \cite{ChenWangWuXu2015,WeiWuXu2026}; angle constraints and the general existence problem are discussed in \cite{MondelloPanov2016,Eremenko2021}. G\'alvez, Hauswirth, and Mira \cite{GHM2013} developed the extrinsic theory in $\E^3$ and showed, in particular, that a peaked sphere with exactly two singularities is rotational. That conclusion concerns closed convex embedded surfaces. Their local theory of immersed singularities is essential for the present treatment of immersions with arbitrary covering multiplicity. Explicit spherical football immersions were constructed in \cite{Wei2025}.

HCMU footballs are basic building blocks in the intrinsic theory of singular extremal K\"ahler metrics on Riemann surfaces. In complex dimension one, Calabi's extremal condition is expressed by the umbilicity of the curvature Hessian \cite{Calabi1982}. The intrinsic theory was developed in \cite{Chen1999,Chen2000,LinZhu2002,WangZhu2000}. Football decompositions and explicit constructions were studied in \cite{ChenChenWu2005,ChenWu2009,ChenWu2011,ChenWuXu2015}. In the extrinsic setting, Peng and Wu \cite{PengWu2020} studied local Weingarten realizations in $\E^3$, and Wang and Wu \cite{WangWu2022} constructed explicit football immersions. Wei and Wu classified local realizations of HCMU metrics as Weingarten surfaces in $\E^3$ \cite{WeiWu2022a} and proved local isometric immersion results in three-dimensional space forms, with an application to $\C P^3$ \cite{WeiWu2022b}.

Our main result is a global rigidity theorem for positive-curvature rotational metrics under an endpoint trace condition. Its key step is to show that this condition forces the mixed coefficient of the shape operator to vanish, without assuming rotational symmetry of the immersion. Reconstruction on the universal cover and the angular closing condition then give a discrete classification. We verify the trace condition for every smooth spherical football immersion. For HCMU footballs, we classify the rotational immersions for the full curvature range and obtain exhaustiveness under the trace condition when the minimum curvature is nonnegative.

Throughout the paper, $S^1=\R/(2\pi\Z)$ and $\N=\{1,2,\ldots\}$. Consider
\begin{equation}\label{eq:metric}
 M=(0,L)\times S^1,\qquad g=\mathrm{d}u^2+f(u)^2\mathrm{d}\theta^2,
\end{equation}
where $0<L<\infty$, $f\in C^\infty(0,L)\cap C^1[0,L]$, $f>0$ on $(0,L)$, and
\begin{equation}\label{eq:endpoints}
 f(0)=f(L)=0,\qquad f'(0)=\alpha>0,\qquad f'(L)=-\beta<0.
\end{equation}
The endpoint derivatives are one-sided. These assumptions yield cone-angle asymptotics with parameters $\alpha$ and $\beta$; the stronger conformal conical regularity is verified below for the footballs under consideration. In the orthonormal frame $e_1=\partial_u$, $e_2=f^{-1}\partial_\theta$, write the shape operator as $S=\left(\begin{smallmatrix}a&b\\b&c\end{smallmatrix}\right)$. The following trace condition is independent of the choice of unit normal.

\begin{theorem}\label{thm:rotational}
Let \eqref{eq:metric}--\eqref{eq:endpoints} define a rotational metric with Gauss curvature $K=-f''/f>0$ on $(0,L)$. Suppose that a smooth isometric immersion into $\E^3$ satisfies
\begin{equation}\label{eq:trace}
 \lim_{u\to0^+}\sup_\theta |f^2b|=0,
 \qquad \lim_{u\to L^-}\sup_\theta |f^2b|=0.
\end{equation}
Then it is congruent to
\begin{equation}\label{eq:model}
 X_m(u,\theta)=\left(\frac{f(u)}m\cos(m\theta),\,
 \frac{f(u)}m\sin(m\theta),\,
 \int_{u_*}^{u}\sqrt{1-\frac{f'(v)^2}{m^2}}\dd v\right),
\end{equation}
where $u_*\in(0,L)$ is fixed and the integer $m\in\N$ is uniquely determined, with
\begin{equation}\label{eq:degree}
 m\geq\max\{\alpha,\beta\}.
\end{equation}
Conversely, every integer satisfying \eqref{eq:degree} defines a smooth isometric immersion of the regular cylinder by \eqref{eq:model}. Each model satisfies \eqref{eq:trace}, since its mixed shape coefficient vanishes identically.
\end{theorem}

The first application gives a classification of all smooth isometric immersions of spherical footballs. Let $g_\alpha$ be a spherical football metric with cone angles $2\pi\alpha$ at both tips. On its regular cylinder,
\begin{equation}\label{eq:sphericalmetric}
 M_\alpha=(0,\pi)\times S^1,\qquad
 g_\alpha=\mathrm{d}u^2+\alpha^2\sin^2u\,\mathrm{d}\theta^2.
\end{equation}
This polar-coordinate expression is justified in Lemma~\ref{lem:footballpolar}.

\begin{theorem}\label{thm:spherical}
Let $\alpha>0$. Every smooth isometric immersion $X:(M_\alpha,g_\alpha)\to\E^3$ is congruent, by a Euclidean rigid motion, to
\begin{equation}\label{eq:sphericalmodel}
 X_{B,\lambda}(u,\theta)=\left(B\sin u\cos(\lambda\theta),\,
 B\sin u\sin(\lambda\theta),\,
 \int_0^u\sqrt{1-B^2\cos^2v}\dd v\right),
\end{equation}
where
\begin{equation}\label{eq:sphericaldegree}
 \lambda\in\N,\qquad \lambda\geq\alpha,\qquad B=\frac\alpha\lambda\in(0,1].
\end{equation}
Conversely, every choice in \eqref{eq:sphericaldegree} defines a smooth isometric immersion of $M_\alpha$. Distinct values of $\lambda$ give noncongruent parametrized immersions.
\end{theorem}

The second application concerns nonconstant HCMU footballs. Orient the meridian so that $K(0)=K_0>K_1=K(L)$, and retain $\alpha=f'(0)$ and $\beta=-f'(L)$. The existence of these endpoint curvature limits follows from the HCMU equation in Section~\ref{sec:hcmu}. A \emph{rotational immersion} is an immersion equivariant under $\theta$-translations and a circle action by Euclidean rotations about a fixed axis. This convention permits nontrivial angular covering degree.

\begin{theorem}\label{thm:hcmu}
Let $g$ be a nonconstant HCMU football of the form \eqref{eq:metric}--\eqref{eq:endpoints}, oriented as above. Then $K_0>0$, $-K_0/2<K_1<K_0$, and $\alpha>\beta>0$.
\begin{enumerate}
\item Its rotational isometric immersions are, up to Euclidean congruence, precisely the maps \eqref{eq:model} with $m\in\N$ and $m\geq\alpha$. In particular, a rotational immersion exists for every HCMU football, including those whose curvature changes sign.
\item A rotational embedding of the regular cylinder exists if and only if $\alpha\leq1$. When it exists, it has $m=1$ and is unique within the rotational class.
\item If $K_1\geq0$, every smooth isometric immersion satisfying \eqref{eq:trace} is one of these rotational models. In particular, this conclusion holds for every immersion with a $C^2$ polar extension at both tips, in the sense of Definition~\ref{def:polar}.
\end{enumerate}
For $K_1<0$, parts \textup{(i)} and \textup{(ii)} classify only the rotational immersions.
\end{theorem}

Section~\ref{sec:prelim} introduces the intrinsic definitions, endpoint conventions, and Gauss--Codazzi equations, and relates the model formulas to earlier rotational constructions. Section~\ref{sec:rigidity} proves Theorem~\ref{thm:rotational} and gives sufficient conditions for the endpoint traces to vanish. Section~\ref{sec:spherical} proves Theorem~\ref{thm:spherical} and describes the geometry of its models. Section~\ref{sec:hcmu} proves Theorem~\ref{thm:hcmu} and discusses endpoint regularity, examples, and local realizations.

\section{Preliminaries}\label{sec:prelim}
\subsection{Conical singularities and spherical footballs}
\begin{definition}[\cite{Troyanov1989}]\label{def:cone}
Let $M$ be a Riemann surface and $p\in M$. A conformal metric $\mathrm{d}s^2$ has a conical singularity at $p$ with cone angle $2\pi\alpha$, where $\alpha>0$ and $\alpha\neq1$, if in a local complex coordinate $z$ centered at $p$ it has the form
\[
 \mathrm{d}s^2=|z|^{2(\alpha-1)}e^{2\psi(z)}|\mathrm{d}z|^2,
\]
where $\psi$ extends continuously to $z=0$.
\end{definition}

Troyanov \cite{Troyanov1989} classified spherical metrics with two conical singularities; another proof is given in \cite{WeiWuXu2026}. We recall the explicit classification because a two-cone metric with integer cone-angle parameter need not be a football in the distance-$\pi$ sense.

\begin{theorem}[\cite{Troyanov1989,WeiWuXu2026}]\label{thm:intrinsic}
Identify $S^2$ with $\C\cup\{\infty\}$. Let $\mathrm{d}s^2$ be a spherical conical metric with singularities at $0$ and $\infty$, both of cone angle $2\pi\alpha$, where $\alpha>0$ and $\alpha\neq1$. Up to a change of coordinate $z\mapsto pz$, $p\in\C\setminus\{0\}$, it has one of the following forms:
\begin{enumerate}[label=\textup{(\arabic*)}]
\item If $\alpha\notin\N$, then
\[
 \mathrm{d}s^2=\frac{4\alpha^2|z|^{2(\alpha-1)}}{(1+|z|^{2\alpha})^2}|\mathrm{d}z|^2.
\]
\item If $\alpha\in\N$, then
\[
 \mathrm{d}s^2=\frac{4\alpha^2|z|^{2(\alpha-1)}}{(1+|z^\alpha+\eta|^2)^2}|\mathrm{d}z|^2,
 \qquad \eta\in\R.
\]
\end{enumerate}
\end{theorem}

\begin{definition}[\cite{WeiWuXu2026}]\label{def:football}
A metric of constant curvature one on a two-sphere with two conical singularities is a \emph{spherical football metric} if the distance between the singularities is $\pi$. The resulting metric sphere is called a \emph{spherical football}. We include the round sphere with a marked antipodal pair as the case $\alpha=1$.
\end{definition}

\begin{lemma}\label{lem:footballpolar}
A spherical football metric with cone angles $2\pi\alpha$ at both tips has the geodesic polar form
\[
 \mathrm{d}s^2=\mathrm{d}r^2+(\alpha\sin r)^2\mathrm{d}\theta^2,
 \qquad (r,\theta)\in(0,\pi)\times S^1.
\]
\end{lemma}
\begin{proof}
For noninteger $\alpha$, the first case of Theorem~\ref{thm:intrinsic} applies. For $\alpha=n\geq2$ an integer, put $F(z)=z^n+\eta$. The metric in the second case is the pullback of the round metric $4|\mathrm{d}w|^2/(1+|w|^2)^2$ by $F$. Every path from $0$ to $\infty$ projects to a path from $\eta$ to $\infty$ of the same length. Conversely, a shortest round-sphere arc joining these two branch values lifts to an arc joining their unique preimages: away from the endpoints, $F$ is an unbranched covering. Hence
\[
 d_{\mathrm{d}s^2}(0,\infty)=d_{\mathrm{round}}(\eta,\infty)
 =\pi-2\arctan|\eta|.
\]
The distance-$\pi$ condition therefore forces $\eta=0$. The case $\alpha=1$ is the round sphere with antipodal markings. Thus the radial formula in the first case applies in every case under consideration.

Set $z=R^{1/\alpha}e^{i\theta}$ and $R=\tan(r/2)$. Then
\[
 |\mathrm{d}z|^2=\frac1{\alpha^2}R^{2/\alpha-2}\mathrm{d}R^2
 +R^{2/\alpha}\mathrm{d}\theta^2,\qquad |z|^{2\alpha}=R^2.
\]
Substitution gives
\[
 \mathrm{d}s^2=\frac4{(1+R^2)^2}\mathrm{d}R^2
 +\frac{4\alpha^2R^2}{(1+R^2)^2}\mathrm{d}\theta^2
 =\mathrm{d}r^2+(\alpha\sin r)^2\mathrm{d}\theta^2.
\]
The radial curves have unit speed, and $r$ is their distance from the first tip.
\end{proof}

\subsection{HCMU footballs}
\begin{definition}[\cite{Chen2000}]\label{def:hcmu}
Let $M$ be a compact Riemann surface, let $p_1,\ldots,p_N\in M$, and set $M^*=M\setminus\{p_1,\ldots,p_N\}$. A conformal metric $g$ on $M^*$ is an \emph{HCMU metric} if its Gauss curvature satisfies
\begin{equation}\label{eq:hcmu}
 K_{,zz}=0,\qquad K\not\equiv\mathrm{constant},
\end{equation}
where the comma denotes covariant differentiation with respect to $g$.
\end{definition}

For $g=e^{2\phi}|\mathrm{d}z|^2$, the covariant derivative is $K_{,zz}=K_{zz}-2\phi_zK_z$. Equivalently,
\[
 \nabla^2K=\tfrac12(\Delta_gK)g,
\]
where $\Delta_g$ is the trace of the Hessian. In this paper, an \emph{HCMU football} is a rotational metric with two tips satisfying this nonconstant-curvature condition. Smooth marked tips are allowed when a cone-angle parameter equals one. In the rotational frame, the Hessian eigenvalues are $K''$ and $(f'/f)K'$, so the HCMU equation becomes $K''=(f'/f)K'$. Its global consequences are derived in Section~\ref{sec:hcmu}.

\subsection{Endpoint and immersion conventions}\label{subsec:conventions}
The metric completion of \eqref{eq:metric}--\eqref{eq:endpoints} adds two points, denoted by $p_0$ and $p_L$. Their intrinsic distances from $(u,\theta)$ are $u$ and $L-u$, respectively. The endpoint slopes give
\[
 f(u)\sim\alpha u\quad(u\to0^+),\qquad
 f(u)\sim\beta(L-u)\quad(u\to L^-).
\]
The area element and Gauss curvature are
\begin{equation}\label{eq:areaK}
 \mathrm{d}A=f\dd u\dd\theta,\qquad K=-\frac{f''}{f}.
\end{equation}
An isometric immersion is a smooth map on the regular cylinder whose first fundamental form is $g$. Neither injectivity nor endpoint regularity is assumed. Two parametrized immersions are \emph{congruent} if they differ by a Euclidean rigid motion, possibly including a reflection. All area and curvature integrals are taken over the parameter domain and therefore count multiplicity.

\begin{definition}[]\label{def:polar}
An immersion has a $C^j$ \emph{polar extension} at a tip if, in the original geodesic polar coordinates $(u,\theta)$, it extends as a $C^j$ map to the corresponding boundary circle of $[0,L]\times S^1$ and is constant on that circle.
\end{definition}
This is a regularity condition on a cylinder with boundary. It does not assert regular immersion through the completed point in a smooth Cartesian coordinate. It also differs from smooth extension in coordinates conformal for the second fundamental form, which is used in Section~\ref{sec:spherical}.

\subsection{The Gauss--Codazzi equations}
Let $X:M\to\E^3$ be an isometric immersion of \eqref{eq:metric}. Choose a global unit normal $N$ and use the convention $S=-\mathrm{d}N$. In the orthonormal frame
\[
 e_1=\partial_u,\qquad e_2=f^{-1}\partial_\theta,\qquad
 S=\begin{pmatrix}a&b\\b&c\end{pmatrix},\qquad H=\tfrac12(a+c),
\]
the connection satisfies
\[
 \nabla_{e_1}e_1=\nabla_{e_1}e_2=0,\qquad
 \nabla_{e_2}e_1=\frac{f'}f e_2,\qquad
 \nabla_{e_2}e_2=-\frac{f'}f e_1.
\]
The Gauss--Codazzi equations are therefore
\begin{align}
 ac-b^2&=K(u),\label{eq:gauss}\\
 b_u&=\frac{a_\theta}{f}-2\frac{f'}f b,\label{eq:codazzi1}\\
 c_u&=\frac{b_\theta}{f}+\frac{f'}f(a-c).\label{eq:codazzi2}
\end{align}
When $K>0$, the eigenvalues of $S$ have the same nonzero sign. Connectedness permits a global choice of normal for which $S>0$. Thus $a,c>0$ and
\begin{equation}\label{eq:bH}
 |b|\leq H,
\end{equation}
because $b^2<ac\leq((a+c)/2)^2$. The fundamental theorem of surface theory determines an immersion on a simply connected domain from its first and second fundamental forms, up to congruence \cite{doCarmo1976,Spivak1999}. For the cylinder, we apply this theorem on the universal cover and then impose the deck-transformation condition.

\subsection{Relation with classical rotational constructions}\label{subsec:classical}
Engman's profile formula \cite[Section~2]{Engman2004} is the case $m=1$ of \eqref{eq:model}. To compare with Giblin--Hwang \cite[equation~(4.1)]{GiblinHwang2004}, introduce $r$ and $\varphi$ by
\[
 \mathrm{d}r=f(u)\,\mathrm{d}u,\qquad \varphi(r(u))=f(u)^2.
\]
Then
\[
 g=\varphi(r)^{-1}\mathrm{d}r^2+\varphi(r)\mathrm{d}\theta^2,
 \qquad \frac{\mathrm{d}\varphi}{\mathrm{d}r}=2f'(u).
\]
Their rotational expression with angular parameter $\kappa$ becomes \eqref{eq:model} with $m$ replaced by $\kappa$, up to interchanging the first two ambient coordinates. Their global ambient-space and period conventions differ from ours; here the target is ordinary Euclidean space and the domain angle has fixed period $2\pi$. Theorem~\ref{thm:rotational} establishes that every immersion satisfying \eqref{eq:trace} has this form and that the closing condition forces the parameter to be an integer.

\section{Proof of Theorem~\ref{thm:rotational}}\label{sec:rigidity}
The proof has two stages. We first derive an elliptic equation for the weighted mixed shape coefficient and use the endpoint traces to show that it vanishes. We then integrate the remaining Gauss--Codazzi equations and impose angular periodicity.

\subsection{An elliptic equation and maximum principles}
\begin{lemma}\label{lem:elliptic}
For an immersion with $K=K(u)>0$, choose the normal so that $S>0$ and put $W=f^2b$. Then
\begin{equation}\label{eq:Wdiv}
 \partial_u\left(\frac ca W_u-\frac{2b}{af}W_\theta\right)
 +\frac1{f^2}W_{\theta\theta}+\frac{f'}f W_u=0.
\end{equation}
For the fixed immersion, this is a linear equation in $W$, with no zeroth-order term, and it is uniformly elliptic on every compact subcylinder.
\end{lemma}
\begin{proof}
Equations \eqref{eq:codazzi1}--\eqref{eq:codazzi2} give
\begin{equation}\label{eq:Wfirst}
 W_u=fa_\theta,\qquad (fc)_u=b_\theta+f'a.
\end{equation}
Since $K_\theta=0$, differentiating \eqref{eq:gauss} yields
\[
 ca_\theta+ac_\theta-2bb_\theta=0.
\]
Substitution of $a_\theta=W_u/f$ and $b_\theta=W_\theta/f^2$ gives
\[
 fc_\theta=\frac{2b}{af}W_\theta-\frac ca W_u.
\]
Differentiate the second equation of \eqref{eq:Wfirst} with respect to $\theta$ and use this identity to obtain \eqref{eq:Wdiv}. Its nondivergence form is
\begin{equation}\label{eq:Wnondiv}
\begin{split}
 0={}&\frac ca W_{uu}-\frac{2b}{af}W_{u\theta}
       +\frac1{f^2}W_{\theta\theta}\\
 &+\left[\left(\frac ca\right)_u+\frac{f'}f\right]W_u
       -\left(\frac{2b}{af}\right)_uW_\theta.
\end{split}
\end{equation}
The symmetric principal matrix is
\begin{equation}\label{eq:principal}
 A=\begin{pmatrix}c/a&-b/(af)\\-b/(af)&1/f^2\end{pmatrix},
 \qquad \det A=\frac{ac-b^2}{a^2f^2}=\frac K{a^2f^2}>0.
\end{equation}
Together with $c/a>0$, this proves positive definiteness. Smoothness and compactness give uniform ellipticity on each compact subcylinder. No uniform ellipticity at the tips is required.
\end{proof}

We recall the weak maximum principle in the form needed below; see \cite[Theorem~3.1 and Corollary~3.2, pp.~32--33]{GilbargTrudinger2001}.

\begin{lemma}[Weak maximum principle]\label{lem:weakmaximum}
Let $D\subset\R^n$, $n\geq2$, be a bounded domain, and consider
\[
 Lv=\sum_{i,j=1}^n A^{ij}(x)\partial_{ij}v
     +\sum_{i=1}^n B^i(x)\partial_i v+C(x)v.
\]
Suppose that $A(x)$ is symmetric and positive definite, with smallest eigenvalue $\lambda(x)>0$, and that
\[
 \sup_{x\in D}\frac{|B^i(x)|}{\lambda(x)}<\infty,
 \qquad i=1,\ldots,n.
\]
For $v\in C^2(D)\cap C^0(\overline D)$, the following statements hold.
\begin{enumerate}
\item If $C\equiv0$, then
\[
 Lv\geq0\ \Longrightarrow\ \max_{\overline D}v=\max_{\partial D}v,
 \qquad
 Lv\leq0\ \Longrightarrow\ \min_{\overline D}v=\min_{\partial D}v.
\]
\item If $C\leq0$, then
\[
 Lv\geq0\ \Longrightarrow\ \max_{\overline D}v\leq\max_{\partial D}v^+,
 \qquad
 Lv\leq0\ \Longrightarrow\ \min_{\overline D}v\geq\min_{\partial D}v^-,
\]
where $v^+=\max\{v,0\}$ and $v^-=\min\{v,0\}$. In particular,
\[
 C\leq0,\quad Lv=0
 \quad\Longrightarrow\quad
 \max_{\overline D}|v|=\max_{\partial D}|v|.
\]
\end{enumerate}
All differential inequalities are imposed in $D$.
\end{lemma}

\begin{lemma}[Maximum principle on a compact cylinder]\label{lem:cylindermaximum}
Let $u_-<u_+$, $\Omega=(u_-,u_+)\times S^1$, and $\overline\Omega=[u_-,u_+]\times S^1$. Consider
\[
 Lv=A^{11}v_{uu}+2A^{12}v_{u\theta}+A^{22}v_{\theta\theta}
       +B^1v_u+B^2v_\theta,
\]
whose coefficients are continuous on $\overline\Omega$ and periodic in $\theta$. Assume that $A=(A^{ij})$ is symmetric and uniformly positive definite: there is a constant $\kappa>0$ such that
\[
 \sum_{i,j=1}^2 A^{ij}(u,\theta)\xi_i\xi_j\geq\kappa|\xi|^2
 \quad\text{for all }(u,\theta)\in\overline\Omega,\ \xi\in\R^2.
\]
If $v\in C^2(\Omega)\cap C^0(\overline\Omega)$ and $Lv\geq0$ in $\Omega$, then
\[
 \max_{\overline\Omega}v
 =\max\left\{\max_\theta v(u_-,\theta),\max_\theta v(u_+,\theta)\right\}.
\]
In particular, if $Lv=0$, then
\[
 \max_{\overline\Omega}|v|
 \leq\max\left\{\max_\theta|v(u_-,\theta)|,\max_\theta|v(u_+,\theta)|\right\}.
\]
\end{lemma}
\begin{proof}
We give a direct argument to account for the periodic coordinate. Choose $k>0$ with $k>\|B^1\|_{L^\infty(\Omega)}/\kappa$ and put $h(u)=e^{ku}$. Then
\[
 Lh=e^{ku}(k^2A^{11}+kB^1)
 \geq e^{ku}\bigl(k^2\kappa-k\|B^1\|_{L^\infty(\Omega)}\bigr)>0.
\]
For every $\delta>0$, the function $v_\delta=v+\delta h$ satisfies $Lv_\delta>0$. At an interior maximum, its gradient would vanish and its Hessian would be negative semidefinite. Positive definiteness of $A$ would then give $Lv_\delta\leq0$, a contradiction. Thus the maximum of $v_\delta$ lies on one of the two end circles; periodicity creates no additional boundary. Letting $\delta\to0^+$ proves the first assertion. Applying it to both $v$ and $-v$ proves the absolute-value estimate.
\end{proof}

\subsection{Proof of the classification}
\begin{proof}[\textbf{Proof of Theorem~\ref{thm:rotational}}]
\emph{Step 1: vanishing of the mixed coefficient.}
Fix $0<\varepsilon<L/2$. By Lemma~\ref{lem:elliptic}, equation \eqref{eq:Wnondiv} satisfies the hypotheses of Lemma~\ref{lem:cylindermaximum} on $[\varepsilon,L-\varepsilon]\times S^1$. The ellipticity constant may depend on $\varepsilon$. Hence
\begin{equation}\label{eq:Wmaximum}
 \sup_{[\varepsilon,L-\varepsilon]\times S^1}|W|
 \leq\max\left\{\sup_\theta|W(\varepsilon,\theta)|,
                  \sup_\theta|W(L-\varepsilon,\theta)|\right\}.
\end{equation}
For any fixed interior point, the right-hand side tends to zero by \eqref{eq:trace}. Therefore $W\equiv0$ and $b\equiv0$. Equation \eqref{eq:Wfirst} implies $a_\theta=0$, and $ac=K(u)>0$ then gives $c_\theta=0$.

\emph{Step 2: integration of the diagonal system.}
The second equation of \eqref{eq:Wfirst} reduces to $(fc)'=f'a$. Consequently,
\[
 (f^2c^2)'=2fc(fc)'=2ff'ac=2ff'K=-2f'f''.
\]
There is a constant $\mu>0$ such that
\begin{equation}\label{eq:firstintegral}
 f'^2+f^2c^2=\mu^2.
\end{equation}
Since $f,c>0$, we have $\mu>|f'|$ in the interior. Taking endpoint limits yields $\mu\geq\max\{\alpha,\beta\}$. The shape coefficients are
\begin{equation}\label{eq:shape}
 c=\frac{\sqrt{\mu^2-f'^2}}f,
 \qquad a=\frac{Kf}{\sqrt{\mu^2-f'^2}}.
\end{equation}

\emph{Step 3: reconstruction and angular closing.}
On the universal cover $(0,L)\times\R$, use \eqref{eq:model} with $m$ replaced by $\mu$. Write $r=f/\mu$ and $q=\sqrt{1-f'^2/\mu^2}$. The derivatives are
\[
 X_u=(r'\cos(\mu\theta),r'\sin(\mu\theta),q),\qquad
 X_\theta=(-f\sin(\mu\theta),f\cos(\mu\theta),0).
\]
Thus $|X_u|=1$, $\langle X_u,X_\theta\rangle=0$, and $|X_\theta|=f$. For the normal
\[
 N=(-q\cos(\mu\theta),-q\sin(\mu\theta),r'),
\]
the mixed second fundamental coefficient vanishes, and direct differentiation gives
\begin{equation}\label{eq:modelII}
 \II=\frac{Kf}{\sqrt{\mu^2-f'^2}}\mathrm{d}u^2
       +f\sqrt{\mu^2-f'^2}\,\mathrm{d}\theta^2.
\end{equation}
For example, $\langle X_{uu},N\rangle=-r''/q=Kf/\sqrt{\mu^2-f'^2}$, using $r'^2+q^2=1$. The model and the lifted original immersion therefore have the same first and second fundamental forms. By the fundamental theorem of surface theory, they are congruent on the universal cover.

The original immersion is $2\pi$-periodic in $\theta$, so its congruent model must also be periodic. Since $f/\mu>0$, this is equivalent to $e^{2\pi i\mu}=1$. Hence $\mu=m\in\N$, proving \eqref{eq:degree} and the asserted formula.

\emph{Step 4: converse and uniqueness.}
Because $f''=-Kf<0$, the function $f'$ decreases strictly from $\alpha$ to $-\beta$. Every integer satisfying \eqref{eq:degree} therefore satisfies $|f'|<m$ in the interior, and \eqref{eq:model} is smooth there. The first fundamental form calculation proves isometry, and $b\equiv0$ gives \eqref{eq:trace}. The height is strictly increasing, and each parallel is traversed exactly $m$ times. Finally, at any fixed interior value of $u$, the positive angular normal curvature is
\[
 c(u)=\frac{\sqrt{m^2-f'(u)^2}}{f(u)}.
\]
It is strictly increasing in $m$ and is invariant under congruence after the normal is chosen so that $S>0$. Distinct integers therefore give distinct congruence classes.
\end{proof}

\subsection{Embeddings and sufficient trace conditions}
\begin{corollary}\label{cor:embedding}
Within the class covered by Theorem~\ref{thm:rotational}, an embedding of the regular cylinder exists if and only if $\max\{\alpha,\beta\}\leq1$. When it exists, it is unique up to congruence and has $m=1$.
\end{corollary}
\begin{proof}
For $m>1$, the distinct points $(u,\theta)$ and $(u,\theta+2\pi/m)$ have the same image. For $m=1$, the strictly increasing height determines $u$ continuously on the image, and the nonzero horizontal component then determines the angular coordinate. Thus the map is an embedding.
\end{proof}

\begin{proposition}\label{prop:traceconditions}
A $C^2$ polar extension at a tip implies the corresponding limit in \eqref{eq:trace}. When $K>0$ and the normal is chosen so that $S>0$, the estimate $f^2\sup_\theta H\to0$ is another sufficient condition. In particular, a uniform bound $H=O(u^{-1})$ at the first tip, or $H=O((L-u)^{-1})$ at the second, suffices.
\end{proposition}
\begin{proof}
Since $\II(\partial_u,\partial_\theta)=\langle X_{u\theta},N\rangle=fb$, we have
\[
 |f^2b|=f|\langle X_{u\theta},N\rangle|\leq f|X_{u\theta}|.
\]
A $C^2$ polar extension makes $X_{u\theta}$ uniformly bounded, while $f\to0$. The mean-curvature criterion follows from \eqref{eq:bH}, and the stated growth bounds follow from $f(u)\sim\alpha u$ and $f(u)\sim\beta(L-u)$.
\end{proof}

\begin{remark}\label{rem:modeltrace}
For the explicit models, $b\equiv0$ yields the trace condition without any estimate on mean-curvature growth or higher endpoint regularity. Proposition~\ref{prop:traceconditions} supplies sufficient conditions for a general immersion to satisfy the same hypothesis.
\end{remark}

\section{Proof of Theorem~\ref{thm:spherical}}\label{sec:spherical}
We apply Theorem~\ref{thm:rotational} after establishing its endpoint trace condition for every spherical football immersion. The main estimate is obtained by comparing the first fundamental form with coordinates conformal for $\II$ near a tip. Throughout this section, the normal is chosen so that $\II>0$; in particular, $H=\tfrac12\tr S\geq1$ for a $K=1$ immersion. The conventions of Section~\ref{subsec:conventions} remain in force.

\subsection{Completion and local singularity theory}
\begin{lemma}[]\label{lem:completion}
Every isometric immersion of \eqref{eq:metric}--\eqref{eq:endpoints} extends continuously to the metric completion. If the two tips have images $P_0$ and $P_L$, then
\[
 |X(u,\theta)-P_0|\leq u,\qquad |X(u,\theta)-P_L|\leq L-u.
\]
Both ends have finite area, counted with multiplicity.
\end{lemma}
\begin{proof}
Isometry gives $|X_u|=1$ and $|X_\theta|=f$. For fixed $\theta$ and $0<u<v<L$,
\[
 |X(v,\theta)-X(u,\theta)|\leq\int_u^v|X_u(s,\theta)|\dd s=v-u.
\]
Thus $X(u,\theta)$ has a limit at each endpoint. Taking the shorter angular arc gives
\[
 |X(u,\theta_1)-X(u,\theta_2)|
 \leq f(u)d_{S^1}(\theta_1,\theta_2)\leq\pi f(u),
\]
where $d_{S^1}$ is the distance on the unit circle. Since $f\to0$, all radial limits at a given endpoint coincide. Radial integration to this common limit yields the asserted estimates, uniformly in $\theta$. Finally,
\[
 \Area\{0<u<\varepsilon\}=2\pi\int_0^\varepsilon f(u)\dd u<\infty,
\]
and the same argument applies at the other end. Equivalently, continuous extension follows from the fact that an isometric immersion is distance nonincreasing from intrinsic distance to Euclidean distance.
\end{proof}

We next recall the results of G\'alvez--Hauswirth--Mira \cite{GHM2013} used in the proof. Their \emph{extrinsic conformal structure} is the conformal structure induced by the positive second fundamental form. 

\begin{lemma}[{\cite[equations~(2.2)--(2.7)]{GHM2013}}]\label{lem:IIconformal}
In a positively oriented $\II$-conformal coordinate $z=s+it$ on a $K=1$ immersion, there are positive functions $\nu,\rho$ and a holomorphic function $Q$ such that
\[
 \I=Q\,\mathrm{d}z^2+2\nu|\mathrm{d}z|^2+\overline Q\,\mathrm{d}\bar z^2,
 \qquad \II=2\rho|\mathrm{d}z|^2,
\]
and
\begin{equation}\label{eq:conformalidentities}
\begin{gathered}
 \rho^2=\nu^2-|Q|^2,\qquad \rho=|Q|\sinh\omega,\\
 H=\frac\nu\rho=\coth\omega,\qquad \Delta\omega=-4\rho<0.
\end{gathered}
\end{equation}
Here $\Delta=\partial_s^2+\partial_t^2$, and the identities involving $\omega$ hold where $Q\neq0$. The zeros of $Q$ are the umbilics; at these points $H=1$, and we set $\omega=+\infty$. Moreover,
\begin{equation}\label{eq:Hsystem}
 X_{ss}+X_{tt}=2X_s\times X_t.
\end{equation}
\end{lemma}

\begin{theorem}[{\cite[Proposition~4 and Theorem~5]{GHM2013}}]\label{thm:singularities}
Let $X:D\setminus\{q\}\to\R^3$ be a $K=1$ immersion of a punctured disk that extends continuously to $q$.
\begin{enumerate}
\item Finite area near $q$ is equivalent to finite total mean curvature there:
\[
 \int\mathrm{d}A<\infty\quad\Longleftrightarrow\quad\int H\dd A<\infty.
\]
\item The following conditions are equivalent:
\begin{enumerate}
\item $H$ is bounded near $q$;
\item the extrinsic conformal end is a punctured disk;
\item $(X,N)$ extends smoothly across $q$ in a disk parametrization;
\item the singularity is removable or is a branch point.
\end{enumerate}
If these conditions fail, the extrinsic conformal end is an annulus.
\end{enumerate}
\end{theorem}

\begin{lemma}[{\cite[Lemma~2]{GHM2013}}]\label{lem:boundaryregularity}
Let $\Omega\subset\R^2$ be a bounded domain with smooth boundary. Suppose that
\[
 Y\in C^\infty(\Omega,\R^3)\cap C^0(\overline\Omega,\R^3)
       \cap W^{1,2}(\Omega,\R^3),\qquad
 Y_{ss}+Y_{tt}=2Y_s\times Y_t.
\]
If $Y|_{\partial\Omega}$ is smooth, then $Y\in C^\infty(\overline\Omega,\R^3)$.
\end{lemma}

The next lemma isolates the boundary extension argument used in the proof of \cite[Theorem~7]{GHM2013}.

\begin{lemma}\label{lem:collar}
For a finite-area end in Theorem~\ref{thm:singularities} with unbounded $H$, a smaller end admits $\II$-conformal coordinates
\[
 \mathcal C=(\R/P\Z)\times(0,t_0],\qquad z=s+it,\qquad P,t_0>0,
\]
in which $X$ extends smoothly to the closed cylinder and $X(s,0)=X(q)$ for every $s$.
\end{lemma}
\begin{proof}
The annular alternative in Theorem~\ref{thm:singularities} gives periodic conformal coordinates of this form. Choose the outer circle in the regular part. The continuous extension at the puncture gives the constant boundary value at $t=0$.

By Lemma~\ref{lem:IIconformal},
\[
 \mathrm{d}A=2\rho\dd s\dd t,\qquad |X_s|^2+|X_t|^2=4\nu.
\]
Theorem~\ref{thm:singularities}\textup{(i)} therefore yields
\begin{equation}\label{eq:energy}
 \frac12\int_{\mathcal C}(|X_s|^2+|X_t|^2)\dd s\dd t
 =\int_{\mathcal C}H\dd A<\infty.
\end{equation}
The continuous extension also makes $X$ bounded, so its components lie in $W^{1,2}(\mathcal C)$. On a conformally equivalent planar annulus, equation \eqref{eq:Hsystem} retains its form, and the data on both boundary circles are smooth. Lemma~\ref{lem:boundaryregularity} gives smoothness up to the boundary. Pulling this extension back to the cylinder proves the assertion.
\end{proof}

\subsection{The endpoint trace and the classification}
The preceding results apply to immersed ends and allow angular covering multiplicity. We now estimate $W=f^2b$ directly.

\begin{lemma}\label{lem:sphericaltrace}
Every smooth isometric immersion of $(M_\alpha,g_\alpha)$, with $f(u)=\alpha\sin u$, satisfies
\begin{equation}\label{eq:traceestimate}
 \sup_\theta|f(u)^2b(u,\theta)|=O(f(u))
 \qquad\text{as }u\to0^+\text{ or }u\to\pi^-.
\end{equation}
In particular, both endpoint limits in \eqref{eq:trace} vanish.
\end{lemma}
\begin{proof}
Fix one tip, and set $r=u$ at the first tip and $r=\pi-u$ at the second. Thus $f=\alpha\sin r$. Lemma~\ref{lem:completion} gives a continuous point extension and finite area. If $H$ is bounded near the tip, then \eqref{eq:bH} gives $|W|\leq f^2H=O(f^2)$, which implies \eqref{eq:traceestimate}.

Suppose instead that $H$ is unbounded, and use the periodic $\II$-conformal collar $(s,t)$ from Lemma~\ref{lem:collar}. The collapsed boundary $t=0$ corresponds to $r=0$. Smooth extension gives $|X_t|\leq M$ for some $M>0$. Expressing the first fundamental form in these coordinates yields
\begin{equation}\label{eq:metriccomparison}
 r_t^2+f^2\theta_t^2=\I(\partial_t,\partial_t)=|X_t|^2\leq M^2,
 \qquad r(s,t)\leq Mt,\qquad f\leq\alpha Mt.
\end{equation}
Indeed, the first identity gives $|r_t|\leq M$; integration from $r(s,0)=0$ and the inequality $\sin r\leq r$ give the remaining estimates. Angular derivatives are taken in a local lift of $\theta$.

The holomorphic differential $Q\,\mathrm{d}z^2$ is not identically zero, since otherwise $H\equiv1$. Its zeros are isolated. Choose $t_0$ so that the outer circle contains no zero of $Q$, and put
\[
 m_0=\min_{s\in\R/P\Z}\omega(s,t_0)>0.
\]
For $0<\delta<t_0$ such that the circle $t=\delta$ contains no zero of $Q$, compare $\omega$ with the harmonic function
\[
 h_\delta(t)=m_0\frac{t-\delta}{t_0-\delta}.
\]
On both boundary circles, $\omega-h_\delta\geq0$. Remove small disks around the finitely many umbilics in the compact subcylinder. The same inequality holds on the added boundary circles when their radii are sufficiently small, because $\omega\to+\infty$ at an umbilic. Away from these disks, \eqref{eq:conformalidentities} gives $\Delta(\omega-h_\delta)=-4\rho<0$. The minimum principle, applied in annular coordinates, yields $\omega\geq h_\delta$. Shrink the disks and let a sequence of admissible $\delta$ tend to zero. We obtain
\begin{equation}\label{eq:omegalower}
 \omega(s,t)\geq kt,\qquad k=\frac{m_0}{t_0}>0.
\end{equation}
At umbilics, we retain the convention $\omega=+\infty$ and $H=1$.

For $x>0$, the inequality $\coth x\leq1+x^{-1}$ implies $\coth(kt)\leq C_0/t$ on $0<t\leq t_0$, where $C_0=t_0+k^{-1}$. Combining \eqref{eq:bH}, \eqref{eq:conformalidentities}, \eqref{eq:metriccomparison}, and \eqref{eq:omegalower}, we find
\begin{equation}\label{eq:directtrace}
 \begin{split}
 |W|&\leq f^2H=f^2\coth\omega\leq f^2\coth(kt)\\
     &\leq C_0\frac{f^2}{t}\leq C_0\alpha M f.
 \end{split}
\end{equation}
The constants are uniform on the periodic collar, which contains a punctured neighborhood of the chosen tip. This proves the uniform-in-$\theta$ estimate \eqref{eq:traceestimate} at either end. Since $f\to0$, both traces vanish. Only an upper bound for $|X_t|$ was used, so the argument also applies when $X_t$ vanishes at points of the collapsed boundary.
\end{proof}

\begin{figure}[tbp]
\centering
\begin{tikzpicture}[
 block/.style={draw=black!65,rounded corners=2pt,fill=black!2,align=center,inner sep=6pt,font=\small},
 wide/.style={block,text width=12.5cm},
 half/.style={block,text width=5.9cm},
 flow/.style={-{Stealth[length=1.5mm]},line width=.5pt,draw=black!75}]
\node[wide,minimum height=1.0cm] (start) at (0,0)
 {Smooth isometric immersion of $(M_\alpha,g_\alpha)$\\
 Continuous completion and finite area (Lemma~\ref{lem:completion})};
\node[half,minimum height=2.0cm] (bounded) at (-3.45,-2.15)
 {Bounded mean curvature\\[3pt]
 $|W|\leq f^2H=O(f^2)$\\at the chosen tip};
\node[half,minimum height=2.0cm] (unbounded) at (3.45,-2.15)
 {Unbounded mean curvature\\
 Smooth $\II$-conformal collar (Lemma~\ref{lem:collar})};
\node[wide,minimum height=1.05cm] (trace) at (0,-4.35)
 {Endpoint trace: $\sup_\theta|f^2b|=O(f)\to0$\\
 Apply the estimates at both tips (Lemma~\ref{lem:sphericaltrace})};
\node[wide,minimum height=1.05cm] (classify) at (0,-6.05)
 {Theorem~\ref{thm:rotational}: $X$ is congruent to $X_m$, $m\in\N$, $m\geq\alpha$\\
 Set $\lambda=m$ and $B=\alpha/m$ to obtain \eqref{eq:sphericalmodel}};
\node[half,minimum height=1.0cm] (converse) at (-3.45,-8.05)
 {Converse\\
 Integer $\lambda$ gives periodicity;\\
 direct differentiation gives\\
 $\I=\mathrm{d}u^2+\alpha^2\sin^2u\,\mathrm{d}\theta^2$};
\node[half,minimum height=1.0cm] (unique) at (3.45,-8.05)
 {Uniqueness\\
 The degree $m$ determines $\lambda$;\\
 distinct degrees give noncongruent\\
 parametrized immersions};
\draw[flow] (start.south)--++(0,-.3)-|(bounded.north);
\draw[flow] (start.south)--++(0,-.3)-|(unbounded.north);
\draw[flow] (bounded.south)--++(0,-.3)-|(trace.north);
\draw[flow] (unbounded.south)--++(0,-.3)-|(trace.north);
\draw[flow] (trace.south)--(classify.north);
\draw[flow] (classify.south)--++(0,-.3)-|(converse.north);
\draw[flow] (classify.south)--++(0,-.3)-|(unique.north);
\end{tikzpicture}
\caption{Proof of Theorem~\ref{thm:spherical}, with $W=f^2b$. The unbounded-mean-curvature case combines the metric comparison \eqref{eq:metriccomparison} with the minimum-principle estimate \eqref{eq:omegalower}. All endpoint estimates are uniform on the corresponding collar.}
\label{fig:proof}
\end{figure}
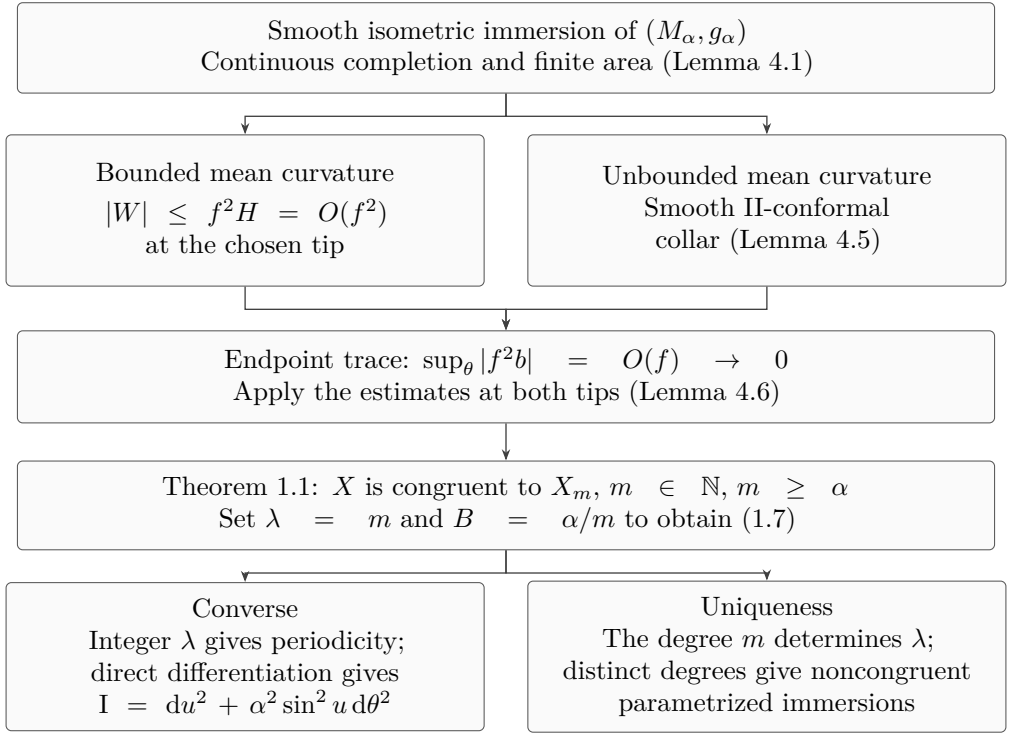

\begin{proof}[\textbf{Proof of Theorem~\ref{thm:spherical}}]
The argument is summarized in Figure~\ref{fig:proof}. Let $X:(M_\alpha,g_\alpha)\to\E^3$ be a smooth isometric immersion. Here $L=\pi$, $f(u)=\alpha\sin u$, $K=1$, and both cone-angle parameters equal $\alpha$. Lemma~\ref{lem:sphericaltrace} verifies the remaining hypothesis of Theorem~\ref{thm:rotational}. Thus $X$ is congruent to \eqref{eq:model} for an integer $m\geq\alpha$. Taking $\lambda=m$ and $B=\alpha/m$ gives \eqref{eq:sphericalmodel}; replacing the lower integration limit by zero changes only the vertical translation.

Conversely, suppose \eqref{eq:sphericaldegree} holds. The integer $\lambda$ makes the formula single valued in $\theta$, and $1-B^2\cos^2u>0$ on $(0,\pi)$. Direct differentiation yields
\[
 \langle X_u,X_u\rangle=1,\qquad
 \langle X_u,X_\theta\rangle=0,\qquad
 \langle X_\theta,X_\theta\rangle=B^2\lambda^2\sin^2u=\alpha^2\sin^2u.
\]
Hence the formula defines a smooth isometric immersion of the regular cylinder. Uniqueness of the degree in Theorem~\ref{thm:rotational} shows that distinct integers $\lambda$ give noncongruent maps.
\end{proof}

\subsection{Geometric consequences}
The classification distinguishes the geometry of the image from the covering multiplicity of its parametrization.

\begin{corollary}\label{cor:sphericalembedding}
A spherical football admits an embedding of its regular cylinder into $\E^3$ if and only if $\alpha\leq1$. Such an embedding is unique up to congruence and is the model with $\lambda=1$.
\end{corollary}
\begin{proof}
By Theorem~\ref{thm:spherical}, every immersion is a model \eqref{eq:sphericalmodel}. For $\lambda>1$, the points $(u,\theta)$ and $(u,\theta+2\pi/\lambda)$ have the same image. For $\lambda=1$, the height is strictly increasing and the horizontal radius is positive, so the inverse on the image is continuous. This gives an embedding precisely when $\lambda=1$ is admissible, that is, when $\alpha\leq1$. Uniqueness follows from the classification.
\end{proof}

\begin{remark}\label{rem:sphericalimage}
The regular image of $X_{B,\lambda}$ is a surface of revolution covered $\lambda$ times by the parameter cylinder. The image has cone-angle parameter $B=\alpha/\lambda$, whereas the intrinsic parameter of the domain is $\alpha$. Indeed, near either tip the image radius is $B\sin r\sim Br$, where $r$ is meridian distance from that tip.
\end{remark}

\begin{remark}\label{rem:smoothtips}
If $B=1$, equivalently $\alpha=\lambda\in\N$, the completed image is the unit sphere, and each regular parallel is traversed $\lambda$ times. If $B<1$, the image has conical tips. In particular, when $\alpha=1$ and $\lambda>1$, the intrinsic metric extends smoothly across the marked tips, but the map does not extend as a smooth regular immersion there. This explains why Theorem~\ref{thm:spherical} is formulated on the regular cylinder.
\end{remark}

\section{Proof of Theorem~\ref{thm:hcmu}}\label{sec:hcmu}
\subsection{The curvature polynomial of an HCMU football}
We derive the curvature polynomial directly from the HCMU equation, fixing the normalization used in the immersion formulas. The resulting description agrees with the intrinsic football constructions in \cite{ChenChenWu2005,ChenWu2011,ChenWuXu2015}.

\begin{proposition}\label{prop:curvaturepolynomial}
Let $g$ be a nonconstant HCMU football of the form \eqref{eq:metric}--\eqref{eq:endpoints}. After reversing $u$ if necessary, there is a constant $\gamma>0$ such that
\begin{equation}\label{eq:Kprime}
 K'=-\gamma f.
\end{equation}
The curvature has finite endpoint limits $K_0=K(0)>K_1=K(L)$. Put $\Sigma=K_0^2+K_0K_1+K_1^2$. Then
\begin{align}
 f'&=\frac{3K^2-\Sigma}{6\gamma},\label{eq:fprime}\\
 f^2&=\frac{(K_0-K)(K-K_1)(K+K_0+K_1)}{3\gamma^2},\label{eq:fpolynomial}\\
 \alpha&=\frac{(K_0-K_1)(2K_0+K_1)}{6\gamma},
 \qquad \beta=\frac{(K_0-K_1)(K_0+2K_1)}{6\gamma}.\label{eq:slopes}
\end{align}
In particular,
\begin{equation}\label{eq:curvaturerange}
 K_0>0,\quad -\frac{K_0}{2}<K_1<K_0,\quad
 \alpha>\beta>0,\quad
 \frac\alpha\beta=\frac{2K_0+K_1}{K_0+2K_1}.
\end{equation}
Conversely, every triple $(K_0,K_1,\gamma)$ satisfying $K_0>0$, $-K_0/2<K_1<K_0$, and $\gamma>0$ determines a football profile with the endpoint slopes in \eqref{eq:slopes}.
\end{proposition}
\begin{proof}
The Hessian equation gives
\[
 K''=\frac{f'}fK',\qquad \left(\frac{K'}f\right)'=0.
\]
Thus $K'/f$ is constant and nonzero, since $K$ is nonconstant. Reversing the meridian coordinate if necessary, write this constant as $-\gamma<0$. Continuity of $f$ on $[0,L]$ and \eqref{eq:Kprime} imply that $K$ has finite endpoint limits. Since $K'<0$ in the interior, these satisfy $K_0>K_1$.

Combining $f''=-Kf$ with \eqref{eq:Kprime} yields
\[
 \left(f'-\frac{K^2}{2\gamma}\right)'=0,
 \qquad f'=\frac{K^2}{2\gamma}+C.
\]
Using $K$ as a coordinate in the interior, we obtain
\[
 \frac{\mathrm{d}(f^2)}{\mathrm{d}K}=-\frac{2f'}\gamma
 =-\frac{K^2}{\gamma^2}-\frac{2C}\gamma.
\]
The endpoint conditions $f(K_0)=f(K_1)=0$ imply
\[
 C=-\frac{K_0^2+K_0K_1+K_1^2}{6\gamma}=-\frac\Sigma{6\gamma}.
\]
Integration, with either endpoint condition, gives \eqref{eq:fpolynomial}, while substitution gives \eqref{eq:fprime}. Evaluation at the endpoints proves \eqref{eq:slopes}. Positivity of $\alpha$ and $\beta$ implies $2K_0+K_1>0$ and $K_0+2K_1>0$. Together with $K_0>K_1$, these inequalities give $K_0>0$ and $K_1>-K_0/2$. Moreover,
\[
 \alpha-\beta=\frac{(K_0-K_1)^2}{6\gamma}>0,
\]
which completes the proof of \eqref{eq:curvaturerange}.

Conversely, suppose the parameters satisfy the stated inequalities, and set
\[
 P(\kappa)=(K_0-\kappa)(\kappa-K_1)(\kappa+K_0+K_1).
\]
Then $P>0$ on $(K_1,K_0)$ and has simple zeros at both endpoints. Define
\begin{equation}\label{eq:quadrature}
\begin{gathered}
 u(K)=\sqrt3\int_K^{K_0}\frac{\mathrm{d}\kappa}{\sqrt{P(\kappa)}},
 \qquad L=\sqrt3\int_{K_1}^{K_0}\frac{\mathrm{d}\kappa}{\sqrt{P(\kappa)}},\\
 f(K)=\frac{\sqrt{P(K)}}{\sqrt3\,\gamma}.
\end{gathered}
\end{equation}
Both endpoint integrals converge. The function $u(K)$ is strictly decreasing and therefore defines $K(u)$ on $(0,L)$. Differentiation gives $K'=-\gamma f$, \eqref{eq:fprime}, and $f''=-Kf$. Hence the metric has Gauss curvature $K$ and satisfies the HCMU equation.

Formula \eqref{eq:fprime} gives exactly the endpoint slopes \eqref{eq:slopes}. The system
\[
 f'=p,\qquad p'=-Kf,\qquad K'=-\gamma f
\]
has an analytic right-hand side. Its finite endpoint values, together with uniqueness for ordinary differential equations, give analytic one-sided extensions. In particular, $f(u)=\alpha u+O(u^3)$ at the first tip, with the analogous expansion at the second. To verify conformal conical regularity, choose a radial conformal coordinate by $\mathrm{d}\log R=\mathrm{d}u/f$. Then
\[
 R=cu^{1/\alpha}(1+O(u^2)),\qquad
 g=(f/R)^2\bigl|\mathrm{d}(Re^{i\theta})\bigr|^2,
\]
which has the form in Definition~\ref{def:cone}. The other endpoint is treated in the same way. This completes the construction; the quadrature also proves uniqueness of the oriented profile.
\end{proof}

\begin{lemma}\label{lem:slopebound}
For the profile in Proposition~\ref{prop:curvaturepolynomial},
\begin{equation}\label{eq:strictslope}
 \sup_{0\leq u\leq L}|f'(u)|=\alpha,
 \qquad |f'(u)|<\alpha\quad(0<u<L).
\end{equation}
These statements remain valid when $K$ changes sign.
\end{lemma}
\begin{proof}
If $K_1\geq0$, then $K>0$ in the interior, so $f''=-Kf<0$. Thus $f'$ decreases strictly from $\alpha$ to $-\beta$, and the assertion follows from $\alpha>\beta$.

Suppose $K_1<0$. As a function of $K$, \eqref{eq:fprime} is a quadratic with minimum $-\Sigma/(6\gamma)$ at $K=0$. Its maximum on $[K_1,K_0]$ is $\alpha$, since $|K_1|<K_0/2<K_0$. Set $t=K_1/K_0\in(-1/2,0)$. Then
\[
 6\gamma\left(\alpha-\frac\Sigma{6\gamma}\right)
 =3K_0^2-2\Sigma=K_0^2(1-2t-2t^2)>0.
\]
Therefore the absolute value of the negative minimum is smaller than $\alpha$. Equality $f'=\alpha$ occurs only at $K=K_0$, that is, at the first tip. This proves both assertions.
\end{proof}

\subsection{Global angular closing and the classification}
\begin{proof}[\textbf{Proof of Theorem~\ref{thm:hcmu}}]
The curvature range and the inequality $\alpha>\beta$ follow from Proposition~\ref{prop:curvaturepolynomial}.

\emph{Part (i): rotational immersions.}
This part requires no sign condition on $K$. A nontrivial circle action by Euclidean rotations about a fixed axis has a nonzero integer angular weight $m$. An ambient reflection reverses its sign, so we may assume $m>0$. After aligning the axis with the third coordinate axis, a general equivariant parametrization has the form
\[
 X(u,\theta)=\bigl(r(u)\cos(m\theta+\psi(u)),\,
 r(u)\sin(m\theta+\psi(u)),\,z(u)\bigr).
\]
Here $r>0$, because $|X_\theta|=f>0$. The angular and mixed metric coefficients give
\[
 m^2r^2=f^2,\qquad mr^2\psi'=0.
\]
Hence $r=f/m$ and $\psi$ is constant; an ambient rotation removes this constant. The meridian coefficient then gives
\begin{equation}\label{eq:zprime}
 z'^2=1-\frac{f'^2}{m^2}.
\end{equation}
In particular, $m\geq\sup|f'|=\alpha$ by Lemma~\ref{lem:slopebound}.

Conversely, if $m\in\N$ and $m\geq\alpha$, the strict interior bound in \eqref{eq:strictslope} makes the right-hand side of \eqref{eq:zprime} positive. Integrating its positive square root gives \eqref{eq:model}, and direct differentiation verifies isometry. For every rotational immersion, $z'$ has fixed sign; reflection in a horizontal plane reduces it to the positive choice, and a vertical translation accounts for the integration constant. Distinct integers yield noncongruent maps because their angular normal curvature magnitudes, $\sqrt{m^2-f'^2}/f$, differ at each fixed interior point. This proves part~\textup{(i)}.

\emph{Part (ii): rotational embeddings.}
For $m>1$, each parallel has multiplicity $m$, so the map is not injective. For $m=1$, the height is strictly increasing and the radius is positive in the interior, which gives an embedding as in Corollary~\ref{cor:embedding}. The degree $m=1$ is admissible precisely when $\alpha\leq1$.

\emph{Part (iii): rigidity under the trace condition.}
If $K_1\geq0$, strict decrease of $K$ gives $K>0$ on $(0,L)$, including the case in which $K$ vanishes only at the final tip. Theorem~\ref{thm:rotational} therefore applies to every immersion satisfying \eqref{eq:trace}, and its degree bound becomes $m\geq\max\{\alpha,\beta\}=\alpha$. A $C^2$ polar extension supplies the trace condition by Proposition~\ref{prop:traceconditions}.

When $K_1<0$, the arguments for parts~\textup{(i)} and~\textup{(ii)} still apply. The elliptic argument does not extend across the zero-curvature parallel: wherever the representation in \eqref{eq:principal} is defined, its determinant vanishes at $K=0$ and is negative where $K<0$. The classification of arbitrary immersions in this range is therefore not a consequence of Theorem~\ref{thm:rotational}.
\end{proof}

\subsection{Endpoint regularity, examples, and local models}
\begin{proposition}\label{prop:endregularity}
Every HCMU model in Theorem~\ref{thm:hcmu} has a smooth one-sided polar extension at both tips. This assertion concerns the parametrization on the cylinder with its boundary circles; it does not assert regular immersion through the completed conical points.
\end{proposition}
\begin{proof}
The analytic profile equations in Proposition~\ref{prop:curvaturepolynomial} give
\begin{equation}\label{eq:endexpansions}
\begin{aligned}
 f(u)&=\alpha u-\frac{K_0\alpha}{6}u^3+O(u^5),
 &K(u)&=K_0-\frac{\gamma\alpha}{2}u^2+O(u^4),\\
 f(L-v)&=\beta v-\frac{K_1\beta}{6}v^3+O(v^5),
 &K(L-v)&=K_1+\frac{\gamma\beta}{2}v^2+O(v^4).
\end{aligned}
\end{equation}
The parity also follows from uniqueness for the analytic system under $u\mapsto-u$ and $f\mapsto-f$ near a tip. If $m>\alpha$, then $1-f'^2/m^2$ is positive at the first tip, so its square root is smooth there. In the borderline case $m=\alpha$, which is possible only when $\alpha$ is an integer,
\[
 1-\frac{f'(u)^2}{\alpha^2}=K_0u^2+O(u^4).
\]
Since $K_0>0$, the positive square root is $u$ times a positive analytic function for $u\geq0$ and therefore has a smooth one-sided extension. At the final tip, $m\geq\alpha>\beta$, so the radicand remains strictly positive. Thus both the radial components and the height integral extend smoothly, while $f=0$ makes each boundary circle map to a single point.
\end{proof}

\begin{example}[Positive curvature]\label{ex:positive}
Take $(K_0,K_1,\gamma)=(2,1,1)$. Then
\[
 f^2=\frac{(2-K)(K-1)(K+3)}3,
 \qquad (\alpha,\beta)=\left(\frac56,\frac23\right).
\]
The $m=1$ model is an embedding of the regular cylinder. Every immersion satisfying the trace condition is rotational, and every embedding satisfying this condition is congruent to the $m=1$ model.
\end{example}

\begin{example}[Changing-sign curvature and multiple coverings]\label{ex:multiple}
Take $(K_0,K_1,\gamma)=(2,-1/2,1)$. Formula~\eqref{eq:slopes} gives
\[
 (\alpha,\beta)=\left(\frac{35}{24},\frac5{12}\right),
 \qquad \min f'=-\frac{13}{24}.
\]
The least admissible degree is $m=2$, for which every parallel is traversed twice; the other rotational immersions have degrees $m=3,4,\ldots$. There is no rotational embedding. Theorem~\ref{thm:hcmu} makes no assertion about nonrotational embeddings in this curvature range.
\end{example}

\begin{example}[An embedded model with changing-sign curvature]\label{ex:embedded}
Keep $K_0=2$ and $K_1=-1/2$, but set $\gamma=35/24$. Then $\alpha=1$ and $\beta=2/7$, so the $m=1$ model is an embedded regular cylinder even though its curvature changes sign.
\end{example}

\begin{proposition}[Analytic local models and the period obstruction]\label{prop:localmodels}
Let $f>0$ be real analytic on an open interval $J$, and let $\mu>0$ satisfy $|f'|<\mu$ on $J$. Formula~\eqref{eq:model}, with $m$ replaced by $\mu$ and a fixed base point $u_*\in J$, defines a real-analytic isometric immersion of $J\times\R$. It descends to $J\times\R/(2\pi\Z)$ if and only if $\mu\in\N$.
\end{proposition}
\begin{proof}
The square root is positive and analytic, and the metric calculation in Section~\ref{sec:rigidity} applies without a sign condition on $K$. The height is independent of $\theta$, while the horizontal component is multiplied by $e^{2\pi i\mu}$ under $\theta\mapsto\theta+2\pi$. Since $f/\mu>0$, descent is equivalent to $e^{2\pi i\mu}=1$, or $\mu\in\N$.
\end{proof}

The HCMU profiles are analytic, so Proposition~\ref{prop:localmodels} applies on every compact interior interval once $\mu$ is sufficiently large. The integer degree is therefore a global closing condition, rather than a consequence of the local immersion equations. The rotational constructions above do not settle the classification of arbitrary immersions when the curvature changes sign. In particular, the present arguments do not determine whether nonrotational immersions exist in this case.

\section*{Declarations}
\noindent\textbf{Generative AI disclosure.} The author used OpenAI's ChatGPT as an assisting tool in developing this paper. The author takes full responsibility for its contents.

\smallskip
\noindent\textbf{Data availability.} This manuscript has no associated data.

\smallskip
\noindent\textbf{Competing interests.} The authors have not disclosed any competing interests.

\end{document}